\documentclass[a4paper,reqno, 11pt]{amsart}
\usepackage[utf8]{inputenc}
\usepackage{tikz}
\usepackage{tikz-cd}
\usepackage[top=2.2cm,bottom=2.2cm,left=2cm,right=2cm]{geometry}
\usepackage{graphicx}
\usepackage{mathtools}

\usepackage{hyperref}

\usepackage[dvips]{epsfig}
\usepackage{color}
\usepackage{verbatim}
\usepackage{amsgen, amstext,amsbsy,amsopn,amssymb, amsthm, mathptmx, amsfonts,amssymb,amscd,amsmath,nicefrac,euscript,enumerate,url,verbatim,calc,framed}
\usepackage[all]{xy}

\theoremstyle{plain}
\newtheorem{thm}{\bf Theorem}[section]

\newtheorem{prop}[thm]{\bf Proposition}
\newtheorem{lem}[thm]{\bf Lemma}

\theoremstyle{definition}

\theoremstyle{remark}
\newtheorem{rem}[thm]{\bf Remark}

\DeclareMathOperator{\coker}{coker}

\DeclareMathOperator{\Sym}{Sym}

\def\height{\operatorname{ht}}

\newsavebox\foobox

\begin{document}

\title{\textbf{Free resolutions of symmetric algebras of a class of ideals of deviation three}}

\author[Neeraj Kumar]{Neeraj Kumar}
\address{Department of Mathematics, Indian Institute of Technology Hyderabad, Kandi, Sangareddy - 502284, INDIA}
\email{neeraj@math.iith.ac.in}

\author[Aniruddha Saha]{Aniruddha Saha}
\address{Department of Mathematics, Indian Institute of Technology Hyderabad, Kandi, Sangareddy - 502284, INDIA}
\email{ma20resch11001@iith.ac.in}

\thanks{\textit{Dates}: Submitted March 5, 2026; Revised June 9, 2026; Accepted June 30, 2026.}
\subjclass[2020]{{Primary 13D02,13C05}; Secondary {13A30} } 
\keywords{Free resolution, Symmetric algebra, Koszul homology}

\begin{abstract}
We study symmetric algebras of graded ideals of deviation three in a polynomial ring. For a class of such ideals, we construct explicit graded free resolutions of the symmetric algebra. Our results extend known explicit constructions, previously available for certain families of ideals of deviation at most two, to a new specific class of ideals of deviation three.
\end{abstract}

\maketitle


\section{Introduction and preliminaries}

The symmetric algebra of an ideal constitutes a fundamental class of bigraded algebras whose graded components correspond to the symmetric powers of the ideal. For a graded ideal, the defining equations of the symmetric algebra can be described explicitly from a presentation matrix of the ideal. Minimal graded free resolutions provide one of the primary tools for studying the algebraic structure and homological invariants of symmetric algebras. However, determining such resolutions for arbitrary ideals remains a challenging problem. Imposing suitable conditions on the generators often leads to deeper structural insight and facilitates a more precise understanding of the associated homological properties.

\vspace{2mm}

Let $A$ be a standard graded polynomial ring in several variables over a field $\mathbf{k}$ of characteristic zero, and let $I\subseteq A$ be a graded ideal. The difference between the minimal number of generators of $I$ and its grade is called the \textit{deviation} of $I$. Ideals of small deviation have been studied extensively, and their symmetric algebras often admit explicit descriptions. In particular, Micali and Roby showed in \cite{MR1971} that the defining equations of the symmetric algebra of an ideal of deviation zero are given by the $2 \times 2$ minors of a $2 \times n$ matrix, where $n$ denotes the minimal number of generators of the ideal. As a consequence, the symmetric algebra in this case admits a graded free resolution given by the Eagon–Northcott complex  \cite{BE1975, EN1962}. 

\vspace{2mm}
	
The approximation complex $\mathcal{Z}_\bullet$ (see Equation (\ref{eq:approximation-complex})) introduced by Herzog, Simis, and Vasconcelos plays a central role in the study of blowup algebras. When $\mathcal{Z}_\bullet$ is acyclic (or under the standard hypotheses ensuring acyclicity), it provides a free resolution of the symmetric algebra of $I$; cf.\ \cite{HSV1,HSV2}. 

\vspace{2mm}

Explicit graded free resolutions of symmetric algebras are also known for graded ideals of deviation one, in both the equigenerated and non-equigenerated settings; see \cite{C19, KV25,RPU25}. For perfect ideals $I \subseteq A$ of deviation two, Cid-Ruiz, Polini, and Ulrich constructed in \cite[Theorem~4.11(ii)]{RPU25} an explicit graded free resolution of $\operatorname{Sym}_A(I)$, under the hypothesis that the associated approximation complex $\mathcal{Z}_\bullet$ is acyclic. Subsequently, Kumar, Saha, and Venugopal \cite{KSV26} extended these results by constructing explicit bigraded free resolutions for certain classes of non-perfect ideals of deviation two.

\vspace{2mm}

Motivated by these developments, in this article we construct explicit graded free resolutions of the symmetric algebras associated with a specific class of ideals of deviation three; see Theorem~\ref{RES}.

\vspace{2mm}

Let $A=\mathbf{k}[x_1,\ldots,x_r]$ be a standard graded polynomial ring, and let $I\subseteq A$ be a graded ideal minimally generated by homogeneous elements $f_1,\ldots,f_n$. Then the symmetric algebra of $I$ is $ \Sym_A(I)=\bigoplus_{j\ge 0}\Sym_A^{j}(I)$. Choose a presentation of $I$ as follows:
\[
A^{t}\xrightarrow{\ \varphi\ } A^{n}\xrightarrow{(f_1,\ldots,f_n)} I\longrightarrow 0,
\]
and write the relation matrix as $\varphi=(a_{ij})$ where $a_{ij}\in A$. Set $S=A[y_1,\ldots,y_n]$. Then $\Sym_A(I)$ admits the presentation $\Sym_A(I)\cong S/J$,
where $J= \langle g_1,\ldots,g_t \rangle\subseteq S$ and
$g_i=\sum_{j=1}^{n} a_{ij}y_j$ for $i=1,\ldots,t$.

\vspace{2mm}

Suppose the deviation of $I$ is three. Let $(\mathbb{F}_\bullet,\varphi_\bullet)$ be a free resolution of $A/I$ of the form
\begin{equation*}
\cdots \longrightarrow A^{\beta_{3}}\stackrel {\varphi_3} \longrightarrow A^{\beta_{2}} \stackrel{\varphi_2}\longrightarrow A^{n} \stackrel{\varphi_1} \longrightarrow A \longrightarrow 0. 
\end{equation*}
Let $(\mathrm{K}_\bullet,\nu_\bullet)$ denote the Koszul complex on the sequence ${f_1,\ldots,f_n}$, equivalently corresponding to the map $\varphi_1$. Then $(\mathrm{K}_\bullet,\nu_\bullet)$ has the form
$$0 \longrightarrow \mathrm{K}_{n} \stackrel{\nu_{n}}\longrightarrow \cdots \stackrel{\nu_{4}}\longrightarrow \mathrm{K}_{3}   \stackrel{\nu_{3}} \longrightarrow \mathrm{K}_{2} \stackrel{\nu_{2}} \longrightarrow \mathrm{K}_{1} \stackrel{\phi_{1}}\longrightarrow A \longrightarrow 0,$$
where $\mathrm{K}_i=\wedge^iA^n$ for $i=1, \ldots, n$. Assume that the approximation complex $\mathcal{Z}_\bullet(I,A)$ is acyclic; then $\mathcal{Z}_\bullet(I,A)$ provides a graded free resolution of the symmetric algebra of $I$, see \cite{HSV1, HSV2}:

\begin{align}
\mathcal{Z}_\bullet(I,A):\quad 
0 &\longrightarrow Z_{n-1}\otimes_A S(-(n-1))
\xrightarrow{\ \partial_{n-1}\ }
\cdots
\xrightarrow{\ \partial_{4}\ }
Z_{3}\otimes_A S(-3)
\xrightarrow{\ \partial_{3}\ }
Z_{2}\otimes_A S(-2) \notag\\
&\xrightarrow{\ \partial_{2}\ }
Z_{1}\otimes_A S(-1)
\xrightarrow{\ \partial_{1}\ }
Z_{0}\otimes_A S
\longrightarrow 0,
\label{eq:approximation-complex}
\end{align}
where $Z_i$ for all $i=1, \ldots, n$, denotes the Koszul cycles of $I$ and the differentials $\partial_{\bullet}$ are induced by the Koszul maps with respect to the sequence $y_1 \ldots, y_n$ \cite{HSV1}. The $i^{\textbf{th}}$ differential map $\partial_{\bullet}^{(i)} : Z_{i} \otimes_A S(-i) \to Z_{i-1} \otimes_A S(-(i-1))$ is given by
     $$ e_{m_1} \wedge \ldots \wedge e_{m_j} \wedge \ldots\wedge e_{m_i} \otimes _A P(\mathbf{y}) \mapsto \sum_{j=1}^i (-1)^{j} g(e_{m_j}) e_{m_1} \wedge \ldots \wedge \hat{e_{m_j}} \wedge \ldots\wedge e_{m_i} \otimes_A y_{j}  P(\mathbf{y}),$$
where $g$ is the corresponding Koszul map of the Koszul complex on the sequence $y_1 \ldots, y_n$. Let $H_i(I,A)$ denote the $i$-th homology module of the Koszul complex. Since the deviation of the ideal $I$ is $3$, the Koszul homology modules satisfy $H_i(I,A) = 0 $ for all $i \geq 4$. Consequently, for each $i \geq 4$, the Koszul complex provides a graded free resolution of the Koszul cycles $Z_i$. This situation is illustrated in the following diagram, where the dotted arrows indicate inclusion maps
\[
\begin{tikzcd}
0 \arrow[r] & \mathrm{K}_{n}  \arrow[r, "\nu_{n}"] & \cdots \arrow[r, "\nu_7"] & \mathrm{K}_6 \arrow[r, "\nu_6"] \arrow[d, two heads] & \mathrm{K}_5 \arrow[r, "\nu_5"] \arrow[d, two heads] & \mathrm{K}_4 \arrow[r, "\nu_3"]  & \cdots \arrow[r] & B \arrow[r] & 0. \\
&  & & Z_5 \arrow[ur, dashrightarrow] & Z_4 \arrow[ur, dashrightarrow] &  &  & & .
\end{tikzcd}
\]
Now consider the following double complex (see \cite[Proposition~4.5]{RPU25}), where $\partial_{\mathbf{y}}$ and $\partial_{\mathbf{f}}$ denote the differentials induced by the Koszul maps with respect to the sequences ${y_1, \ldots, y_n}$ and $f_1,\ldots,f_n$, respectively. Moreover, each column of this double complex is exact.
	

	\begin{equation} \label{double_complex}
		\begin{tikzpicture}[baseline=(current  bounding  box.center)]
			
			\matrix (m) [matrix of math nodes,row sep=2.5em,column sep=2.3em,minimum width=2em, text height=1.8ex, text depth=0.25ex]				{ 
				\vdots & \vdots & \vdots & \vdots  \\
				0 & \wedge^{n}A^{n} \otimes_AS(-(n-2)) & \cdots & \wedge^{6}A^{n} \otimes_AS(-4) \\
				\wedge^{n}A^{n} \otimes_AS(-(n-1)) & \wedge^{n-1}A^{n} \otimes_AS(-(n-2)) & \cdots & \wedge^{5}A^{n} \otimes_AS(-4)  \\
				Z_{n-1} \otimes_A S(-(n-1)) & Z_{n-2} \otimes_A S(-(n-2)) & \cdots & Z_{4} \otimes_A S(-4) \\
				0 & 0 &  &  0  \\
			};						
			\path[-stealth]
			
			
			(m-1-2) edge node [right] {} (m-2-2)
			(m-1-4) edge node [right] {} (m-2-4)
			(m-2-1) edge node [right] {} (m-3-1)
			(m-2-2) edge node [right] {$\partial_\mathbf{f}$} (m-3-2)
			(m-2-4) edge node [right] {$\partial_\mathbf{f}$} (m-3-4)
			(m-3-1) edge node [right] {$\partial_\mathbf{f}$} (m-4-1)
			(m-3-2) edge node [right] {$\partial_\mathbf{f}$} (m-4-2)
			(m-3-4) edge node [right] {$\partial_\mathbf{f}$} (m-4-4)
			(m-4-1) edge node [right] {} (m-5-1)
			(m-4-2) edge node [right] {} (m-5-2)
			(m-4-4) edge node [right] {} (m-5-4)

			(m-2-1) edge node [above] {} (m-2-2)
			(m-2-2) edge node [above] {$\partial_\mathbf{y}$} (m-2-3)
			(m-2-3) edge node [above] {$\partial_\mathbf{y}$} (m-2-4)
			
			(m-3-1) edge node [above] {$\partial_\mathbf{y}$} (m-3-2)
			(m-3-2) edge node [above] {$\partial_\mathbf{y}$} (m-3-3)
			(m-3-3) edge node [above] {$\partial_\mathbf{y}$} (m-3-4)
			
			(m-4-1) edge node [above] {$\partial_\mathbf{y}$} (m-4-2)
			(m-4-2) edge node [above] {$\partial_\mathbf{y}$} (m-4-3)
			(m-4-3) edge node [above] {$\partial_\mathbf{y}$} (m-4-4)	
;		\end{tikzpicture}
	\end{equation}
Let $(\mathbb{C}_\bullet,\sigma_\bullet)$ denote the total complex associated with the truncation of the double complex \eqref{double_complex}. Its modules are of the form 
	\begin{equation*} \label{truncated_complex}
	C_i= \oplus_{j=0}^{i}\wedge^{5+i}A^{n} \otimes_AS(-j-4)
\end{equation*}
for $i = 0, \ldots, n-5$, obtained by taking the direct sums along the diagonals of the truncated double complex. By \cite[Proposition~4.5]{RPU25}, one obtains a complex $\mathcal{L}_\bullet$ whose homology is $\operatorname{Sym}_A(I)$, and whose tail consists of graded free $S$-modules
\begin{equation*} \label{NewComplex}
		\mathcal{L}_\bullet:\quad 0 \longrightarrow C_{n-5} \stackrel{\sigma_{n-5}} \longrightarrow C_{n-4} \stackrel{\sigma_{n-4}}\longrightarrow \cdots \stackrel{\sigma_1} \longrightarrow C_{0}   \stackrel{{\eta}}\longrightarrow Z_{3} \otimes_A S(-3) \longrightarrow Z_{2} \otimes_A S(-2) \stackrel{} \longrightarrow  Z_{1} \otimes_A S(-1) \stackrel{}\longrightarrow S \longrightarrow 0.
\end{equation*}
Moreover, $\mathcal{L}_\bullet$ is acyclic if and only if  the approximation complex  $\mathcal{Z}_\bullet(I,A)$ is acyclic \cite[Proposition~4.5]{RPU25}. Hence, under our assumption that $\mathcal{Z}_\bullet(I,A)$ is acyclic, the complex $\mathcal{L}_\bullet$ provides a resolution of $\operatorname{Sym}_A(I)$. Therefore, once graded free resolutions of 
	\[
	Z_{3}\otimes_A S(-3),\qquad
	Z_{2}\otimes_A S(-2),\qquad \text{ and} \qquad
	Z_{1}\otimes_A S(-1)
	\]
 are known, an explicit graded free resolution of $\operatorname{Sym}_A(I)$ can be obtained by iteratively applying the mapping cone construction.

In this article, we construct a class of ideals of deviation $3$ for which the resolution  $\operatorname{Sym}_A(I)$ can be made explicit; see Section~\ref{Classofideals} and Theorem~\ref{RES}.

\section{Homological properties of Certain Classes of Deviation Three Ideals }

In this section, we develop homological tools for constructing graded free resolutions of the symmetric algebras of a certain class of ideals of deviation $3$. Let $A=\mathbf{k}[x_1,\ldots,x_r]$ be a graded polynomial ring. Let $I^{\prime\prime\prime} \subseteq A $ be an ideal of deviation $3$. We define a class of deviation-three ideals by adjoining, in an iterative manner, a sequence that is regular on $A/I^{\prime\prime\prime}$ as well as on $H_2(I^{\prime\prime\prime},A)$. The following lemmas give conditions under which these homological properties are preserved under this construction.



\begin{lem}\label{exact-H}
Let $I\subseteq A$ be an ideal, and let $h\in A$. Fix an integer $i\ge 0$. Assume that $h$ is $H_i(I,A)$-regular. Set $J=I+ \langle h \rangle$. Then
\[
H_{i+1}(J,A)\cong H_{i+1}(I,A)\big/ h\,H_{i+1}(I,A).
\]
\end{lem}

\begin{proof} Consider the short exact sequence of Koszul complexes
\begin{align*}
0 \longrightarrow \mathrm{K}_\bullet(I, A) \xrightarrow{i} \mathrm{K}_\bullet(J, A) \xrightarrow{\pi} \mathrm{K}_\bullet(I, A)[-1] \longrightarrow 0.
\end{align*}
It induces a long exact sequence in Koszul homology
\begin{align*}
\cdots \longrightarrow 	H_{i+1}(I, A) \xrightarrow{\cdot {h}} H_{i+1}(I, A)  \longrightarrow H_{i+1}(J, A) \longrightarrow H_{i}(I, A) \xrightarrow{\cdot {h}} H_{i}(I, A) \longrightarrow \cdots \longrightarrow H_0(J, A) \longrightarrow 0.
\end{align*}
Since $h$ is $H_i(I,A)$-regular, the map $H_i(I,A)\xrightarrow{\cdot {h}}H_i(I,A)$ is injective. Hence the connecting map $H_{i+1}(J,A)\to H_i(I,A)$ is zero, and we obtain an exact sequence
\begin{align*}
H_{i+1}(I, A) \xrightarrow{\cdot {h}} H_{i+1}(I, A) \longrightarrow H_{i+1}(J, A) \longrightarrow 0.
\end{align*}
Therefore, we have $$ H_{i+1}(J, A) \cong H_{i+1}(I, A)/ h  H_{i+1}(I, A),$$ as claimed.
\end{proof}

\begin{lem}\label{exact-H-1} 
Let $I'''\subseteq A$ be an ideal of deviation $3$, and let $h\in A$. Assume that $h$ is regular on $A/I'''$, on $H_1(I''',A)$, and on $H_2(I''',A)$, and set $J=I'''+\langle h \rangle$. Let $(\mathrm{L}_\bullet,\tau_\bullet)$ be a free resolution of $H_3(I''',A)$. Define a complex $(\mathrm{L}'_\bullet,\tau'_\bullet)$ by
\[
\mathrm{L}'_0=\mathrm{L}_0,\qquad \mathrm{L}'_i=\mathrm{L}_i\oplus \mathrm{L}_{i-1}\ \ (i\ge 1),
\]
together with
\[
\tau'_1=\begin{pmatrix}-\tau_1 & h\,\mathrm{id}_{\mathrm{L}_0}\end{pmatrix},
\qquad
\tau'_{i+1}=
\begin{pmatrix}
-\tau_{i+1} & h\,\mathrm{id}_{\mathrm{L}_i}\\
0 & \tau_i
\end{pmatrix}\ \ (i\ge 1).
\]
Then $(\mathrm{L}'_\bullet,\tau'_\bullet)$ is a free resolution of $H_3(J,A)$.
\end{lem}

\begin{proof} There is a natural short exact sequence of Koszul complexes
\begin{equation}\label{eq:ses_koszul_mapping cone}
0 \longrightarrow \mathrm{K}_\bullet(I''', A) \xrightarrow{i} \mathrm{K}_\bullet(J, A) \xrightarrow{\pi} \mathrm{K}_\bullet(I''', A)[-1] \longrightarrow 0.
\end{equation}
Since both $I'''$ and $J$ have deviation $3$, we have $H_t(I''',A)=0=H_t(J,A)$ for all $t\ge 4$. Hence \eqref{eq:ses_koszul_mapping cone} induces the following long exact sequence in Koszul homology:
\begin{align*}
0 &\longrightarrow H_3(I''',A) \xrightarrow{\cdot h} H_3(I''',A)
\longrightarrow H_3(J,A) \longrightarrow H_2(I''',A) \xrightarrow{\cdot h} H_2(I''',A) \longrightarrow H_2(J,A) \\
&\longrightarrow H_1(I''',A) \xrightarrow{\cdot h} H_1(I''',A)
\longrightarrow H_1(J,A) \longrightarrow H_0(I''',A) \xrightarrow{\cdot h} H_0(I''',A) \longrightarrow H_0(J,A) \longrightarrow 0.
\end{align*}
Since $h$ is regular on $A/I'''=H_0(I''',A)$, on $H_1(I''',A)$, and on $H_2(I''',A)$, Lemma~\ref{exact-H} shows that, for $t=1,2,3$,
\[
H_t(J,A)\cong H_t(I''',A)\big/ h\,H_t(I''',A).
\]
In particular, applying this with $t=3$ gives
\[
H_3(J,A)\cong H_3(I''',A)\big/ h\,H_3(I''',A).
\]
Thus $H_3(J,A)$ is the cokernel of multiplication by $h$ on $H_3(I''',A)$.  Since $(\mathrm{L}_\bullet,\tau_\bullet)$ be a free resolution of $H_3(I''',A)$, then multiplication by $h$ on $H_3(I''',A)$ lifts to a chain map $ \mathrm{L}_\bullet \xrightarrow{\cdot h} \mathrm{L}_\bullet$. Let $(\mathrm{L}'_\bullet,\tau'_\bullet)$ be the mapping cone of this chain map. By construction,
\[
\mathrm{L}'_0=\mathrm{L}_0,\qquad \mathrm{L}'_i=\mathrm{L}_i\oplus \mathrm{L}_{i-1}\ \ (i\ge 1),
\]
and the differentials are given by
\[
\tau'_1=\begin{pmatrix}-\tau_1 & h\,\mathrm{id}_{\mathrm{L}_0}\end{pmatrix},
\qquad
\tau'_{i+1}=
\begin{pmatrix}
-\tau_{i+1} & h\,\mathrm{id}_{\mathrm{L}_i}\\
0 & \tau_i
\end{pmatrix}\ \ (i\ge 1).
\]
Since the mapping cone resolves the cokernel of multiplication by $h$ on $H_3(I''',A)$, it follows that $\mathrm{L}'_\bullet$ is a free resolution of
\[
\coker\!\big(H_3(I''',A)\xrightarrow{\cdot h}H_3(I''',A)\big)\cong H_3(J,A).
\]
\end{proof}

\begin{lem}\label{snake-ses-reg}
    Let $0 \to U \to V \to W \to 0$ be a short exact sequence of $A$-modules. Suppose $f \in A$ is regular on $U$ and $W$. Then $f$ is also regular on $V$, and the induced sequence $0 \to U/fU \to V/fV \to W/fW \to 0$ is exact.
\end{lem}
\begin{proof}
    Now consider the following commutative diagram based on the given short exact sequence of $A$-modules
    \begin{equation}
    \begin{tikzcd}
     0 \arrow[r] & U \arrow[r, "\psi"] \arrow[d, "{{\cdot f}}" left] \arrow[d, "\alpha'' " right]  & V \arrow[r, "\epsilon"] \arrow[d, "{{\cdot f}}" left] \arrow[d, "\beta'' " right] & W \arrow[r] \arrow[d, "{\cdot f}" left] \arrow[d, "\gamma'' " right] & 0 \\
     0 \arrow[r] & U \arrow[r, "\psi"] & V \arrow[r, "\epsilon"] & W \arrow[r] & 0.
   \end{tikzcd}
    \end{equation}
  Applying the Snake Lemma to the commutative diagram of the short exact sequence above, we obtain the following long exact sequence of  $A$-modules
  \begin{align*}
   0 \longrightarrow \ker(\alpha'') \longrightarrow \ker(\beta'') \longrightarrow \ker(\gamma'') \longrightarrow \operatorname{coker}(\alpha'') \longrightarrow \operatorname{coker}(\beta'') \longrightarrow \operatorname{coker}(\gamma'') \longrightarrow 0.
   \end{align*}
   Since $f \in A$ is regular on $U$ and $W$ then $\ker(\alpha'')= 0 =\ker(\gamma'')$. Therefore $\ker(\beta'')=0$, i.e., $f$ is regular on $V$. Moreover we have 
   \[
   0 \longrightarrow \operatorname{coker}(\alpha'') \longrightarrow \operatorname{coker}(\beta'') \longrightarrow \operatorname{coker}(\gamma'') \longrightarrow 0.
   \]    
    Therefore the sequence $0 \to U/fU \to V/fV \to W/fW \to 0$ is exact.    
\end{proof}

Let $ I_0'''=\langle h_1,\ldots,h_n,z_1,z_2,z_3 \rangle \subseteq A $ be an ideal of deviation $3$, where $h_1,\ldots,h_n$ is an $A$-regular sequence. Set
\[
I_0=\langle h_1,\ldots,h_n\rangle,\qquad
I_0'=\langle h_1,\ldots,h_n,z_1\rangle,\qquad
I_0''=\langle h_1,\ldots,h_n,z_1,z_2\rangle.
\]
Assume that the colon ideals satisfy
\[
(I_0:z_1)\subseteq (I_0:z_2)\subseteq (I_0:z_3).
\]
We now define, inductively, a class of ideals $\{I_k'''\}_{k\ge 1}$ obtained from $I_0'''$ by adjoining elements of $A$.

\vspace{2mm}

\begin{description}\label{Classofideals}
\item[\textbf{Class I}]\label{ClassI}
For $k\ge 1$, set
\[
I_k'''=I_0'''+ \langle f_1,\ldots,f_k \rangle,
\]
where $f_1,\ldots,f_k\in A$. Moreover, the sequence $f_1,\ldots,f_k$ satisfies the following conditions:
\begin{enumerate}
\item $f_1,\ldots,f_k$ is $A/I_0$-regular.
\item $f_1,\ldots,f_k$ is $A/I_0'$-regular.
\item $f_1,\ldots,f_k$ is $A/I_0''$-regular.
\item \label{C-4} $f_1,\ldots,f_k$ is $A/I_0'''$-regular.
\item\label{4} $f_1,\ldots,f_k$ is $M$-regular, where $M=(I_0':z_2)/ I_0'$.
\item\label{5} $f_1,\ldots,f_k$ is $M'$-regular, where $M'=(I_0'':z_3)/I_0''$.
\item\label{6} $f_1,\ldots,f_k$ is $N$-regular, where $N=H_1(I_0''',A)$.
\end{enumerate}
\end{description}

\begin{prop}\label{H1-H2-regular}
Let $\{I_k'''\}_{k\ge 1}$ be a family of ideals of \rm{\textbf{Class~I}}. Then the sequence $f_1,\ldots,f_k$ is simultaneously regular on $H_1(I_0'',A)$ and $H_2(I_0''',A)$.
\end{prop}

\begin{proof} It has already been noted that $I_0'''$ is an ideal of deviation $3$ in $A$, and that
\[
I_0=\langle h_1,\ldots,h_n\rangle,\qquad 
I_0'=\langle h_1,\ldots,h_n,z_1\rangle,\qquad 
I_0''=\langle h_1,\ldots,h_n,z_1,z_2\rangle
\]
are subideals of $I_0'''$ of deviations $0$, $1$, and $2$, respectively. We first consider the short exact sequence of Koszul complexes associated to the inclusion $I_0\subseteq I_0'$:
\[
0 \longrightarrow \mathrm{K}_\bullet(I_0, A) \xrightarrow{i} \mathrm{K}_\bullet(I_0', A) \xrightarrow{\pi} \mathrm{K}_\bullet(I_0, A)[-1] \longrightarrow 0.
\]
Since $I_0$ and $I_0'$ have deviations $0$ and $1$, respectively, the induced long exact sequence in Koszul homology yields
\begin{equation}\label{Kos-dev-1}
0 \longrightarrow H_1(I_0',A) \longrightarrow H_0(I_0,A) \xrightarrow{{\cdot z_1}} H_0(I_0,A)
\longrightarrow H_0(I_0',A) \longrightarrow 0.
\end{equation}
Here we use that $H_1\!\big(K_\bullet(I_0,A)[-1]\big)=H_0(I_0,A)\cong A/I_0$. 
By exactness of \eqref{Kos-dev-1}, we obtain
\begin{equation} \label{LES1}
H_1(I_0', A) \cong \ker({\cdot{z_{1}}}) = (0 :_{A/I_0} z_{1}).
\end{equation}
Next, consider the short exact sequence of Koszul complexes associated to ideals $I_0'$ and $I_0''$:
\begin{align}\label{Kos-dev-1'}
0 \longrightarrow \mathrm{K}_\bullet(I_0', A) \xrightarrow{i} \mathrm{K}_\bullet(I_0'', A) \xrightarrow{\pi} \mathrm{K}_\bullet(I_0', A)[-1] \longrightarrow 0.
\end{align}
The corresponding long exact sequence in Koszul homology contains the segment
\begin{align}\label{LES-H1}
0 \longrightarrow  H_2(I_0'', A) \longrightarrow H_1(I_0', A) \xrightarrow{\cdot{z_2}} H_1(I_0', A)  \xrightarrow{\psi} H_1(I_0'', A) \xrightarrow{\epsilon}  \\ 
H_0(I_0', A) \xrightarrow{\cdot{z_2}} H_0(I_0', A) \longrightarrow H_0(I_0'', A) \longrightarrow 0.
\end{align}
From \eqref{LES1} and the above long exact sequence \eqref{LES-H1}, we obtain
\begin{equation}\label{LES2}
H_2(I_0'', A) \cong \ker({\cdot{z_2}}) = (0 :_{A/I_0} z_{1}) \cap (0 :_{A/I_0} z_2).
\end{equation}
Since $(I_0:z_1)\subseteq (I_0:z_2)$, it follows that
\begin{equation}\label{H2-I_0''}
    H_2(I_0'',A) \cong (0:_{A/I_0}z_1).
\end{equation}
Moreover, the map $H_1(I_0',A)\xrightarrow{{\cdot z_2}}H_1(I_0',A)$ of the long exact sequence \eqref{LES-H1} is identically zero, so that
\[
\ker(\psi)=\operatorname{im}(\cdot z_2)=0.
\]
Moreover, $\operatorname{im}(\epsilon) = \ker({\cdot{z_2}}) = (0 :_{A/I_0'} z_2)$. Therefore from the previous long exact sequence \eqref{LES-H1} one can obtain the following short exact sequence,
\begin{align}\label{ses1}
0 \longrightarrow H_1(I_0', A) \xrightarrow{\psi} H_1(I_0'', A) \xrightarrow{\epsilon} (0 :_{A/I_0'} z_2) \longrightarrow 0.
\end{align}
We next consider the following commutative diagram of short exact sequences.

\begin{equation}\label{CD1}
\begin{tikzcd}
0 \arrow[r] 
& H_1(I_0', A) 
  \arrow[r, "\psi"] 
  \arrow[d, "{{\cdot f_1}}" left]
  \arrow[d, " \alpha " right]
& H_1(I_0'', A) 
  \arrow[r, "\epsilon"] 
  \arrow[d, "{{\cdot f_1}}" left]
  \arrow[d, "{\beta}" right]
& (0 :_{A/I_0'} z_2) 
  \arrow[r] 
  \arrow[d, "{{\cdot f_1}}" left]
  \arrow[d, "{\gamma}" right]
& 0 \\
0 \arrow[r] 
& H_1(I_0', A) 
  \arrow[r, "\psi"] 
& H_1(I_0'', A) 
  \arrow[r, "\epsilon"] 
& (0 :_{A/I_0'} z_2) 
  \arrow[r] 
& 0
\end{tikzcd}
\end{equation}
Applying the Snake Lemma to the commutative diagram of short exact sequences above, we obtain the following long exact sequence of  $A$-modules
\begin{align}\label{snake-lemma}
0 \longrightarrow \ker(\alpha) \longrightarrow \ker(\beta) \longrightarrow \ker(\gamma) \longrightarrow \operatorname{coker}(\alpha) \longrightarrow \operatorname{coker}(\beta) \longrightarrow \operatorname{coker}(\gamma) \longrightarrow 0.
\end{align}
Note that
\begin{equation}\label{end}
\ker(\alpha) = (0 :_{A/I_0} z_{1}) \cap (0 :_{A/I_0} f_1) \quad \mbox{[by \eqref{LES1}]}.
\end{equation} and
\begin{equation*}
\ker(\gamma) = (0 :_{A/I_0'} z_2) \cap (0 :_{A/I_0'} f_1)=0 \quad \mbox{[as $f_1$ is $A/I_0'$-regular]}.
\end{equation*}
From the hypothesis, the element $f_1$ is regular over $A/I_0$ and $A/I_0'$. Therefore, injectivity of the maps $\alpha$ and $\gamma$ and Lemma~\ref{snake-ses-reg}, implies that  $\beta$ is also injective.  In other words, $f_1$ is $H_1(I_0'',A)$-regular. Thus, from the long exact sequence \eqref{snake-lemma} and using Lemma~\ref{snake-ses-reg} we obtain the following exact sequence:
\begin{align}\label{ses-coker}
0 \longrightarrow \frac{H_1(I_0', A)}{f_1 H_1(I_0', A)} \longrightarrow \frac{H_1(I_0'', A)}{f_1 H_1(I_0'', A)} \longrightarrow \frac{(0 :_{A/I_0'} z_2)}{f_1 (0 :_{A/I_0'} z_2)} \longrightarrow 0.
\end{align}
Since $f_1$ is regular on $A/I_0'$ and $A/I_0''$, by Lemma~\ref{exact-H}, the short exact sequence \eqref{ses-coker} becomes the following:
\begin{align*}
0 \longrightarrow H_1((I_0',f_1), A) \longrightarrow H_1((I_0'',f_1), A) \longrightarrow \frac{(0 :_{A/I_0'} z_2)}{f_1 (0 :_{A/I_0'} z_2)} \longrightarrow 0.
\end{align*}
Here the ideal $(I_0',f_1)=(I_0,z_1,f_1)=(I_0,f_1,z_1)$ also have deviation $1$. By using similar argument one can observe that $$H_1(I_0',A)/ f_1 H_1(I_0',A) \cong H_1((I_0',f_1), A) \cong (0 :_{A/(I_0, f_1)} z_{1}).$$ Since $f_2$ is regular over $A/(I_0,f_1)$, then $f_2$ is regular over $H_1(I_0',A)/ f_1 H_1(I_0',A)$. Therefore by using the similar argument repeatedly one obtains $f_2, \ldots, f_k$ is regular on $H_1(I_0',A)/ f_1 H_1(I_0',A)$. Since $f_2, \ldots, f_k$ is also regular on ${(0 :_{A/I_0'} z_2)}/{f_1 (0 :_{A/I_0'} z_2)}$, using Lemma~\ref{snake-ses-reg} repeatedly, it follows that $f_2, \ldots, f_k$ is regular on $H_1(I_0'',A)/ f_1  H_1(I_0'',A)$. Thus $f_1, \ldots, f_k$ is regular on $H_1(I_0'',A)$.
Next, consider the short exact sequence of Koszul complexes associated to the ideals $I_0''$ and $I_0'''$:
\begin{align*}
	0 \longrightarrow \mathrm{K}_\bullet(I_0'', A) \xrightarrow{i} \mathrm{K}_\bullet(I_0''', A) \xrightarrow{\pi} \mathrm{K}_\bullet(I_0'', A)[-1] \longrightarrow 0.
\end{align*}
This induces the following long exact sequence in Koszul homology,
\begin{align*}
   0 \longrightarrow  H_3(I_0''', A) \longrightarrow H_2(I_0'', A) \xrightarrow{\cdot{z_3}} H_2(I_0'', A)  \xrightarrow{\psi_2}  H_2(I_0''', A) \xrightarrow{\epsilon_2}  H_1(I_0'', A) \xrightarrow{\cdot{z_3}} \\ H_1(I_0'', A)  \xrightarrow{\psi_1} H_1(I_0''', A) \xrightarrow{\epsilon_1}  H_0(I_0'', A) \xrightarrow{\cdot{z_3}} H_0(I_0'', A) \longrightarrow H_0(I_0''', A) \longrightarrow 0.
\end{align*}
From this sequence and Equation $(\ref{LES2})$, we obtain, 
\begin{equation}\label{LES3}
H_3(I_0''', A)  \cong \ker({\cdot{z_3}}) = (0 :_{A/I_0} z_{1}) \cap (0 :_{A/I_0} z_2) \cap (0 :_{A/I_0} z_3).
\end{equation}
Since $(I_0: z_1 ) \subseteq (I_0: z_2) \subseteq (I_0: z_3)$, it follows that
$$
H_3(I_0''', A) \cong H_2(I_0'', A) \cong (0 :_{A/I_0} z_{1}), \quad \text{and} \quad  \ker(\psi_2) = \operatorname{im}({\cdot{z_3}}) = 0.
$$
Now consider the following short exact sequences:
\begin{equation}\label{s-e-1}
 0 \to H_2(I_0'', A)  \xrightarrow{\psi_2}  H_2(I_0''', A) \xrightarrow{\epsilon_2} \operatorname{im} (\epsilon_2) \to 0,
\end{equation}

\begin{equation}\label{s-e-2}
 0 \to \operatorname{im} (\epsilon_2) \to  H_1(I_0'', A) \xrightarrow{\cdot{z_3}} \operatorname{im} ({\cdot{z_3}}) \to 0,
\end{equation}

\begin{equation}\label{s-e-3}
0 \to \operatorname{im} ({\cdot{z_3}}) \to H_1(I_0'', A)  \xrightarrow{\psi_1} \operatorname{im} (\psi_1) \to 0,
\end{equation}

\begin{equation}\label{s-e-4}
0 \to  \operatorname{im} (\psi_1) \to H_1(I_0''', A) \xrightarrow{\epsilon_1} \operatorname{ker} {(\cdot z_3)} \to 0.
\end{equation}
Moreover, $\operatorname{im}(\epsilon_1) = \ker({\cdot{z_3}}) = (0 :_{A/I_0''} z_3)$. Therefore, from the short exact sequence \eqref{s-e-4} we obtain the following commutative diagram:
\begin{equation}\label{CD2}
\begin{tikzcd}
0 \arrow[r] 
& \operatorname{im}(\psi_1) 
  \arrow[r] 
  \arrow[d, "{{\cdot f_1}}" left]
  & H_1(I_0''',A) 
  \arrow[r, "\epsilon_1"] 
  \arrow[d, "{{\cdot f_1}}" left]
  & (0 :_{A/I_0''} z_3) 
  \arrow[r] 
  \arrow[d, "{{\cdot f_1}}" left]
  & 0 \\
0 \arrow[r] 
& \operatorname{im}(\psi_1) 
  \arrow[r] 
& H_1(I_0''',A) 
  \arrow[r, "\epsilon_1"] 
& (0 :_{A/I_0''} z_3) 
  \arrow[r] 
& 0.
\end{tikzcd}
\end{equation}
Since $f_1$ is $H_1(I_0''',A)$-regular, it is also regular on $\operatorname{im}(\psi_1)$.
Moreover, $f_1$ is regular on $(0:_{A/I_0''}z_3)$. Applying the Snake Lemma to the above commutative diagram we obtain  the following short exact sequence:
$$
0 \longrightarrow \frac{\operatorname{im} (\psi_1)}{f_1 \operatorname{im} (\psi_1)} \longrightarrow \frac{H_1(I_0''', A)}{f_1 H_1(I_0''', A)} \longrightarrow \frac{(0 :_{A/I_0''} z_3)}{f_1 (0 :_{A/I_0''} z_3)} \longrightarrow 0.
$$
Since $f_1,\ldots,f_k$ is regular on $H_1(I_0''',A)$ and on $(0:_{A/I_0''}z_3)$, it follows that $f_1,\ldots,f_k$ is regular on $\operatorname{im}(\psi_1)$. Using the same argument as above, and since $f_1,\ldots,f_k$ is regular on $H_1(I_0'',A)$, we further deduce that $f_1,\ldots,f_k$ is regular on $\operatorname{im}(\cdot z_3)$ and on $\operatorname{im}(\epsilon_2)$, by \eqref{s-e-2} and \eqref{s-e-3}, respectively. Finally, from the short exact sequence \eqref{s-e-1} we obtain the following commutative diagram:
	
\begin{equation*}\label{CD2}
\begin{tikzcd}
 0 \arrow[r] 
 & H_2(I_0'', A) 
 \arrow[r, "\psi_2"] 
 \arrow[d, "{{\cdot f_1}}"left]
 \arrow[d, "\alpha' " right] 
 & H_2(I_0''', A) 
 \arrow[r, "\epsilon_2"] 
 \arrow[d, "{{\cdot f_1}}"left] 
 \arrow[d, "\beta' " right]
 & \operatorname{im} (\epsilon_2) 
 \arrow[r] \arrow[d, "{\cdot f_1}"left] 
 \arrow[d, "\gamma' " right]
 & 0 \\
0 \arrow[r] 
& H_2(I_0'', A) 
\arrow[r, "\psi_2"] 
& H_2(I_0''', A) 
\arrow[r, "\epsilon_2"]  
& \operatorname{im} (\epsilon_2) 
\arrow[r] 
& 0.
\end{tikzcd}
\end{equation*}
From the given hypothesis the element $f_1$ is regular on $A/I_0$ and on $A/I_0'$. Hence by \eqref{H2-I_0''} it is clear that the map $\alpha'$ is injective. Moreover, injectivity of the maps $\alpha'$ and $\gamma'$ and using Lemma~\ref{snake-ses-reg}, one obtains  $\beta'$ is also injective. In other words, $f_1$ is $H_2(I_0''',A)$-regular. Using Lemma~\ref{snake-ses-reg} we obtain the following short exact sequence:
\begin{align}\label{ses-coker_1}
0 \longrightarrow \frac{H_2(I_0'', A)}{f_1 H_2(I_0'', A)} \longrightarrow \frac{H_2(I_0''', A)}{f_1 H_2(I_0''', A)} \longrightarrow \frac{\operatorname{im} (\epsilon_2) }{f_1 \operatorname{im} (\epsilon_2) } \longrightarrow 0.
\end{align}
Since $f_1$ is regular on $A/I_0'$ and $A/I_0''$, by Lemma~\ref{exact-H}, the short exact sequence \eqref{ses-coker_1} becomes the following:
\begin{align*}
0 \longrightarrow H_2((I_0'',f_1), A) \longrightarrow H_2((I_0''',f_1), A) \longrightarrow \frac{\operatorname{im} (\epsilon_2)}{f_1 \operatorname{im} (\epsilon_2)} \longrightarrow 0.
\end{align*}
Here the ideal $(I_0'',f_1)=(I_0,z_1,z_2,f_1)=(I_0,z_1,f_1,z_2)=(I_0,f_1,z_1,z_2)$ also have deviation $2$. Using similar arguments, one can observe that $$H_2(I_0'',A)/ f_1 H_2(I_0'',A) \cong H_2((I_0'',f_1), A) \cong (0 :_{A/(I_0, f_1)} z_{1}) \cap (0 :_{A/(I_0,f_1)} z_2).$$ Since $f_2$ is regular over $A/(I_0,f_1)$, then $f_2$ is regular over $H_2(I_0'',A)/ f_1 H_2(I_0'',A)$. Therefore by using the similar arguments repeatedly one obtains $f_2, \ldots, f_k$ is regular on $H_2(I_0'',A)/ f_1 H_2(I_0'',A)$. Since $f_2, \ldots, f_k$ is also regular on ${\operatorname{im} (\epsilon_2)}/{f_1 \operatorname{im} (\epsilon_2)}$, using Lemma~\ref{snake-ses-reg} repeatedly, it follows that $f_2, \ldots, f_k$ is regular on $H_2(I_0''',A)/ f_1  H_2(I_0''',A)$. Thus $f_1, \ldots, f_k$ is regular on $H_2(I_0''',A)$.
\end{proof}

\begin{rem}
    In the proof of Proposition~\ref{H1-H2-regular} the condition \eqref{C-4} of the \textbf{Class I} has not been used.
\end{rem}

\section{ Free Resolutions of symmetric algebra of Deviation $3$ ideals}

Free resolutions of symmetric algebras are known in several cases: deviation $0$ (complete intersections)  \cite{BE1975, EN1962}, deviation $1$ \cite{C19,RPU25,KV25}, and perfect ideals of deviation $2$ \cite{RPU25}; see also \cite{KSV26} for classes of non-perfect deviation $2$ ideals. In this section, we initiate a corresponding study for non-perfect ideals of deviation $3$. Our main result below constructs graded free resolutions of $\Sym_A(I_k''')$ for a family $\{I_k'''\}_{k\ge 0}$ of deviation $3$ ideals described in \textbf{Class I}.

\subsection*{Set-up}
Let $A=\mathbf{k}[X]$ be a polynomial ring, where $X:= x_1, \ldots, x_r$. Let $\{I_k'''\}_{k\ge 0}$ be a family of ideals as in \textbf{Class~I}. We fix the following graded free complexes over $A$.
\begin{enumerate}[\rm (a)]
\item A minimal graded free resolution of $I_0'''$,
\begin{equation}\label{res_I0}
\mathcal{F}^{(0)}_{\bullet}:\quad \cdots 
\xrightarrow{\varphi^{(0)}_3} A^{\beta_{2}} 
\xrightarrow{\varphi^{(0)}_2} A^{\beta_{1}} 
\xrightarrow{\varphi^{(0)}_1} A \longrightarrow 0,
\end{equation}
where $\beta_i$ are the Betti numbers of $A/I_0'''$.

\item A minimal graded free resolution of the third Koszul homology module $H_3(I_0''',A)$,
\begin{equation}\label{res_H2}
\mathcal{L}^{(0)}_{\bullet}:\quad \cdots 
\xrightarrow{\tau^{(0)}_3} A^{\beta'_{3}} 
\xrightarrow{\tau^{(0)}_2} A^{\beta'_{2}}  
\xrightarrow{\tau^{(0)}_1} A^{\beta'_{1}} \longrightarrow 0,
\end{equation}
where $\beta'_i$ are the Betti numbers of $H_3(I_0''',A)$.

\item The Koszul complex on $I_k'''$,
\begin{equation}\label{res_Koszul}
\mathrm{K}^{(k)}_{\bullet}:\quad \cdots 
\xrightarrow{\nu^{(k)}_3} A^{\beta''_{2}} 
\xrightarrow{\nu^{(k)}_2} A^{\beta''_{1}}  
\xrightarrow{\nu^{(k)}_1} A \longrightarrow 0.
\end{equation}

\item For each $n\ge 0$, a graded free resolution of the second Koszul cycle module $Z^{(n)}_2$ of $I_n'''$ in the family of ideals $\{I_k'''\}_{k \geq 0}$,
\begin{equation}\label{res_Z2}
\mathfrak{C}^{(n)}_{\bullet}:\quad \cdots 
\xrightarrow{\vartheta^{(n)}_3} A^{\beta'''_{3}} 
\xrightarrow{\vartheta^{(n)}_2} A^{\beta'''_{2}}  
\xrightarrow{\vartheta^{(n)}_1} A^{\beta'''_{1}} \longrightarrow 0.
\end{equation}
\end{enumerate}

\begin{thm}\label{RES}
Let $A=\mathbf{k}[X]$ and let $\{I_k'''\}_{k\ge 0}$ be a family of ideals as in \rm{\textbf{Class I}}. Assume that the approximation complex $\mathcal{Z}_{\bullet}(I_k''',A)$ is acyclic. Let $S=\mathbf{k}[X,Y]$, where $Y=[y_1,\ldots,y_{\mu(I_k''')}]$.  Then the symmetric algebra $\Sym_A(I_k''')$ admits a graded free resolution over $S$ that can be written explicitly in terms of the complexes fixed in the above set-up, namely  \eqref{res_I0}, \eqref{res_H2}, \eqref{res_Koszul}, and \eqref{res_Z2}.
\end{thm}

\begin{proof} 
Consider the $\mathcal{Z}$-complex associated to the ideal $I_k'''$, namely the approximation complex
\begin{align*}
\mathcal{Z}_\bullet(I_k''',A): \hspace{2 mm} 0 \longrightarrow Z_{n-1}(I_k''',A) \otimes S(-(n-1))  \stackrel{\partial_{n-1}^{(k)}}\longrightarrow \cdots \stackrel{\partial_{4}^{(k)}}\longrightarrow Z_{3}(I_k''',A) \otimes S(-3)  \stackrel{\partial_{3}^{(k)}} \longrightarrow Z_{2}(I_k''',A) \otimes S(-2) \stackrel{\partial_{2}^{(k)}} \longrightarrow  \\
Z_{1}(I_k''',A) \otimes S(-1)\stackrel{\partial_{1}^{(k)}}\longrightarrow  Z_0(I_k''',A) \otimes S \longrightarrow 0
\end{align*}
where $Z_i(I_k''',A)$, for $i=1,\ldots,n$, denotes the $i$th Koszul cycle module of $I_k'''$. Since the approximation complex $\mathcal{Z}_\bullet(I_k''',A)$ is acyclic, it follows that $\mathcal{Z}_\bullet(I_k''',A)$ is a graded $S$-resolution of the symmetric algebra $\Sym(I_k''')$
\begin{equation} \label{NewComplex}
0 \longrightarrow C_{n-5}^{(k)} \stackrel{\sigma_{n-5}^{(k)}} \longrightarrow  \cdots \stackrel{\sigma_1^{(k)}} \longrightarrow C_{0}^{(k)}   \stackrel{\eta^{(k)}} \longrightarrow Z_{3}(I_k''',A) \otimes_A S(-3) \stackrel{\partial_3^{(k)}} \longrightarrow Z_{2}(I_k''',A) \otimes_A S(-2) \stackrel{\partial_2^{(k)}} \longrightarrow  Z_{1}(I_k''',A) \otimes_A S(-1) \stackrel{\partial_1^{(k)}}\longrightarrow S \longrightarrow 0.
\end{equation}
Next, consider the short exact sequence obtained from the acyclic complex \eqref{NewComplex},
	\begin{equation} \label{SES1}
	0 \longrightarrow \operatorname{coker}(\sigma_1^{(k)})   \stackrel{\eta^{(k)}} \longrightarrow Z_{3}(I_k''',A) \otimes_A S(-3) \stackrel{\partial_3^{(k)}} \longrightarrow \operatorname{im}(\partial_3^{(k)})
	\longrightarrow 0.
\end{equation}
Using Proposition~\ref{H1-H2-regular} and Lemma~\ref{exact-H}, we obtain, for each $0\le j<k$,
\[
H_3(I_{j+1}''',A)\cong H_3(I_j''',A)\big/ f_{j+1}\,H_3(I_j''',A).
\]
Since $f_1, \ldots, f_k$ is regular on $A/I_0'''$ and $H_1(I_0''', A)$, applying the mapping cone construction (as shown in Lemma~\ref{exact-H-1}) to the short exact sequences 
\[ 0 \longrightarrow H_3(I_j''', A) \xrightarrow{\cdot f_{j+1}} H_3(I_j''', A) \longrightarrow H_3(I_{j+1}''', A) \longrightarrow 0 \]
iteratively, we obtain a graded free resolution of $H_3(I_k''',A)$ in terms of the free resolution \eqref{res_H2} of $H_3(I_0''',A)$. Consequently, we obtain a graded free resolution
\[
\mathcal{L}^{(k)}_\bullet:\quad \cdots \longrightarrow \mathcal{L}^{(k)}_2
\longrightarrow \mathcal{L}^{(k)}_1 \longrightarrow \mathcal{L}^{(k)}_0 \longrightarrow 0
\]
of $H_3(I_k''',A)$, where for each $\ell\ge 0$ the free module $\mathcal{L}^{(k)}_\ell$ can be written explicitly in terms of the free modules appearing in the resolution of $H_3(I_0''',A)$. We now consider another short exact sequence
$$
0 \longrightarrow \operatorname{im(\nu^{(k)}_4)} \longrightarrow Z_3(I_k''', A) \longrightarrow H_3(I_k''', A) \longrightarrow 0,
$$
and observe that the truncated Koszul complex
$$
0 \longrightarrow \mathrm{K}_{n}^{(k)}\stackrel{\nu^{(k)}_n}  \longrightarrow \cdots \stackrel{\nu^{(k)}_5} \longrightarrow \mathrm{K}_{4}^{(k)} \longrightarrow 0,
$$
is a free resolution of $\operatorname{im}\!\big(\nu^{(k)}_4\big)$. Applying the Horseshoe Lemma to the above short exact sequence, we obtain a free resolution
$\mathcal{G}^{(k)}_\bullet$ of $Z_3(I_k''',A)$ over $A$ such that, for each $\ell\ge 0$,
\[
\mathcal{G}^{(k)}_\ell = K^{(k)}_{\ell+4}\oplus \mathcal{L}^{(k)}_\ell.
\]
Tensoring with $S(-3)$ provides a graded free resolution
\[
\widetilde{\mathcal{G}}^{(k)}_\bullet := \mathcal{G}^{(k)}_\bullet \otimes_A S(-3)
\]
of $Z_3(I_k''',A)\otimes_A S(-3)$ over $S$. Now, taking the mapping cone of the map $\eta^{(k)}$ (Equation \eqref{SES1}) produces a graded free resolution of $\operatorname{im}(\partial_3^{(k)})$:
\[ (\text{Cone}(\eta^{(k)}))_\ell = C_{\ell-1}^{(k)} \oplus (\mathrm{K}_{\ell+4}^{(k)} \oplus \mathcal{L}_\ell^{(k)}) \otimes S(-3) \]
Next, consider the short exact sequence
\begin{equation}\label{SES2}
0 \longrightarrow \operatorname{im}(\partial_3^{(k)}) \xrightarrow{\eta'^{(k)}}
Z_{2}(I_k''',A)\otimes_A S(-2) \xrightarrow{\partial_2^{(k)}}
\operatorname{im}(\partial_2^{(k)}) \longrightarrow 0,
\end{equation}
where $\eta'^{(k)}$ denotes the inclusion map. Tensoring \eqref{res_Z2} with $S(-2)$ we obtain a graded free resolution of $Z_2(I_k''', A) \otimes S(-2)$ over $S$, namely
\[
\widetilde{\mathfrak{C}}^{(k)}_\bullet := \mathfrak{C}^{(k)}_\bullet \otimes_A S(-2),
\qquad\text{so that}\qquad
\widetilde{\mathfrak{C}}^{(k)}_\ell = \mathfrak{C}^{(k)}_\ell \otimes_A S(-2)\ \ \text{for all }\ell\ge 0.
\]
Taking the mapping cone of $\eta'^{(k)}$ in \eqref{SES2} gives a graded free resolution of $\operatorname{im}(\partial_2^{(k)})$ over $S$. In particular,
$$(\text{Cone}(\eta'^{(k)}))_\ell=\text{Cone}(\eta^{(k)})_{\ell-1} \oplus \tilde{\mathfrak{C}}^{(k)}_\ell= C_{\ell-2}^{(k)} \oplus (\mathrm{K}_{\ell+3}^{(k)} \oplus \mathcal{L}_{\ell-1}^{(k)}) \otimes S(-3) \oplus {\mathfrak{C}}^{(k)}_{\ell} \otimes S(-2).$$
We next consider the short exact sequence
\begin{equation}\label{SES3}
0 \longrightarrow \operatorname{im}(\partial_2^{(k)}) \xrightarrow{\eta''^{(k)}}
Z_{1}(I_k''',A)\otimes_A S(-1) \xrightarrow{\partial_1^{(k)}}
\operatorname{im}(\partial_1^{(k)}) \longrightarrow 0,
\end{equation}
where $\eta''^{(k)}$ denotes the inclusion map. Since the sequence $f_1,\ldots,f_k$ is $A/I_0'''$-regular, we may apply the mapping cone construction iteratively to obtain a graded free resolution of $I_k'''$ in terms of the given graded free resolution of $I_0'''$.
Moreover, since $Z_1(I_k''',A)$ is isomorphic to the first syzygy module of $I_k'''$, a truncation of the resulting resolution yields a graded free resolution of $Z_1(I_k''',A)$ of the form
\[
\cdots \longrightarrow F^{(k)}_4 \xrightarrow{\varphi^{(k)}_4} F^{(k)}_3 
\xrightarrow{\varphi^{(k)}_3} F^{(k)}_2 \longrightarrow 0.
\]
Let $\mathcal{H}^{(k)}_\bullet$ denote a graded free resolution of $Z_1(I_k''',A)$, where
\[
\mathcal{H}^{(k)}_\ell = F^{(k)}_{\ell+2}\qquad \text{for all }\ell\ge 0.
\]
Let $\widetilde{\mathcal{H}}^{(k)}_\bullet$ be the induced graded free resolution of
$Z_1(I_k''',A)\otimes_A S(-1)$, obtained by tensoring with $S(-1)$, i.e.,
\[
\widetilde{\mathcal{H}}^{(k)}_\ell = F^{(k)}_{\ell+2}\otimes_A S(-1)\qquad \text{for all }\ell\ge 0.
\]
Taking the mapping cone of $\eta''^{(k)}$ in Equation \eqref{SES3}, we obtain a graded free resolution of the defining ideal of the symmetric algebra of $I_k'''$. More precisely,
$$
(\text{Cone}(\eta''^{(k)}))_\ell=\text{Cone}(\eta'^{(k)})_{\ell-1} \oplus \tilde{\mathcal{H}}_{\ell}^{(k)}  = C_{\ell-3}^{(k)} \oplus (\mathrm{K}_{\ell+2}^{(k)} \oplus \mathcal{L}_{\ell-2}^{(k)}) \otimes S(-3) \oplus {\mathfrak{C}}^{(k)}_{\ell-1} \otimes S(-2) \oplus {F}_{\ell+2}^{(k)}  \otimes S(-1).
$$ 
for all $\ell\ge 0$.
\end{proof}

We close this section by emphasizing the nature of the input required in Theorem~\ref{RES}.
First, one must work with a family of deviation $3$ ideals $\{I_k'''\}_{k\ge 0}$ obtained from the inductive construction.
Second, the resulting graded free $S$-resolution of $\Sym(I_k''')$ is expressed in terms of the following data:
a minimal graded free resolution of $I_0'''$ \eqref{res_I0}, a minimal graded free resolution of the third Koszul homology module $H_3(I_0''',A)$ \eqref{res_H2}, and, for each $n\ge 0$, a graded free resolution of the second Koszul cycle module $Z_2^{(n)}$ of $I_n'''$ \eqref{res_Z2}.
In addition, the inductive process depends on Proposition~\ref{H1-H2-regular}, which ensures that $f_1,\ldots,f_k$ is simultaneously regular on $H_1(I_0'',A)$ and on $H_2(I_0''',A)$.

While these hypotheses are natural from a homological viewpoint, they are not always convenient in applications, since the relevant Koszul homology and cycle modules are typically difficult to control explicitly, even for fairly concrete ideals.
Consequently, direct applications of Theorem~\ref{RES} are not immediate in general. Accordingly, Theorem~\ref{RES} should be viewed as a structural result: it provides a conceptual and, in principle, explicit description of a graded free resolution of $\Sym(I_k''')$ once the relevant Koszul homology and cycle modules are understood.

In the next subsection, we impose additional hypotheses that allow some of these modules to be described more concretely, leading to an explicit graded free resolution of the symmetric algebra of a specific deviation $3$ ideal.

\subsection{Application}

In this subsection, we study perfect ideals of deviation $3$ in a regular local ring $A$. Let $J\subseteq A$ be such an ideal. It is well known that
\[
\omega_{A/J}\cong H_3(J,A),
\]
and therefore a free resolution of $H_3(J,A)$ is determined by a free resolution of $J$.  Using this observation, we construct an explicit free resolution of the symmetric algebra $ \Sym(J, h)$, where $J$ is a height $2$ perfect ideal of deviation $3$ in $A$ and $h$ is $A/J$-regular. Assume that $J$ is minimally generated by $\underline{t}=\{t_1,t_2,z_1,z_2,z_3\}$.  
Set
\[
J_0=(t_1,t_2),\qquad J_1=(t_1,t_2,z_1),\qquad J_2=(t_1,t_2,z_1,z_2).
\]
Here, $J_0$ is generated by an $A$-regular sequence, $J_1$ has deviation $1$, and $J_2$ has deviation $2$. 
Since $J$ is a height $2$ perfect ideal, the Hilbert--Burch theorem shows that the minimal free resolution of $A/J$ has the form
\begin{align*}
0 \longrightarrow A^{4}
\xlongrightarrow{\varphi}
A^{5}
\xlongrightarrow{[\underline{t}]}
A
\longrightarrow 0.
\end{align*}
Note that $h$ is $A/J$-regular. Hence $\langle J,h \rangle$ is a perfect ideal of height $3$. Consider the short exact sequence
\[
0 \longrightarrow A/J \xrightarrow{\cdot h} A/J \longrightarrow A/\langle J,h \rangle \longrightarrow 0 .
\]
Taking the mapping cone of this multiplication map ${\cdot h}$ provides a free resolution of $A/ \langle J,h \rangle$:
\begin{align*}
	0 \longrightarrow A^{4} 
	\longrightarrow A^{4} \oplus A^{5} 
	\longrightarrow A^{5} \oplus A 
	\longrightarrow A \longrightarrow 0.
\end{align*}
By \cite[Corollary~3.3.9]{BH98}, a free resolution of the canonical module $\omega_{A/(J+(h))}$ is obtained from a free resolution of $A/(J+(h))$; in particular, it has the following form:
$$
0 \longrightarrow A \longrightarrow A \oplus A^{5} \longrightarrow A^{5} \oplus A^{4} \longrightarrow A^{4} \longrightarrow 0. 
$$
Since $\langle J,h \rangle$ is a perfect ideal in the regular local ring $A$ of deviation $3$, we have
\[
\omega_{A/{\langle J,h \rangle}} \cong H_3({\langle J,h \rangle}, A).
\]
Moreover, $\mu(\langle J,h \rangle)=6$ and $\height (\langle J,h \rangle) = 3$. Hence, by \cite[Theorem~1.6.16]{BH98},
\[
H_i(\langle J,h \rangle,A)=0 \qquad \text{for all } i=4,5,6.
\]
Now consider the short exact sequence
$$
0 \longrightarrow B_3(\langle J,h \rangle, A) \longrightarrow Z_3(\langle J,h \rangle, A) \longrightarrow H_3(\langle J,h \rangle, A) \longrightarrow 0,
$$
where $B_3(\langle J,h \rangle, A)= \operatorname{im}(\nu_4)$. Since 
the deviation of $\langle J,h \rangle$ is $3$, the truncated Koszul complex of $\langle J,h \rangle$,
\[
0 \longrightarrow A \stackrel{\nu_6}  \longrightarrow A^6 \stackrel{\nu_5} \longrightarrow A^{15} \longrightarrow 0,
\]
is a free resolution of $B_3(\langle J,h \rangle, A)$. Therefore, applying the Horseshoe Lemma to the above short exact sequence yields a free resolution of
$Z_3(\langle J,h \rangle,,A)$, given by the complex
$$
0 \longrightarrow A \longrightarrow A \oplus(A \oplus A^{5}) \longrightarrow A^6 \oplus(A^{5} \oplus A^{4}) \longrightarrow A^{15} \oplus A^{4} \longrightarrow 0. 
$$
Since $\mu(\langle J,h \rangle)=6$ by \cite[Proposition~4.5]{RPU25}, the symmetric algebra $\Sym(\langle J,h \rangle)$ is resolved by the following complex, whose tail consists of graded free $S$-modules:
\begin{equation*} 
\mathcal{L}_\bullet:\quad 0 \longrightarrow C_{1} \stackrel{\sigma_1} \longrightarrow C_{0}   \stackrel{}\longrightarrow \tilde{Z_{3}} \otimes_A S(-3) \stackrel{i} \longrightarrow \tilde{Z_{2}} \otimes_A S(-2) \stackrel{j} \longrightarrow  \tilde{Z_{1}} \otimes_A S(-1) \stackrel{k} \longrightarrow S \longrightarrow 0
\end{equation*}
Here we set $\tilde{Z_{i}}= Z_i(\langle J,h \rangle, A)$ for $i=1,2,3$. Moreover,
\[
C_0=\bigwedge\nolimits^{5} A^{6}\otimes_A S(-4)\cong A^{6}\otimes_A S(-4)\cong S(-4)^6,
\]
and
\[
C_1=\bigwedge\nolimits^{6} A^{6}\otimes_A S(-4)\ \oplus\ \bigwedge\nolimits^{6} A^{6}\otimes_A S(-5)
\cong S(-4)\oplus S(-5).
\]
Assume that the approximation complex $\mathcal{Z}_\bullet(\langle J,h \rangle,A)$ is acyclic. Then, by \cite[Proposition~4.5]{RPU25}, the complex $\mathcal{L}_\bullet$ is also acyclic and provides a graded $S$-free resolution of the symmetric algebra of the ideal $\langle J,h \rangle$. We now consider the short exact sequences induced by $\mathcal{L}_\bullet$:
\begin{equation}\label{first-seq}
0 \to \operatorname{coker}(\sigma_1) \stackrel{}\longrightarrow \tilde{Z_{3}} \otimes_A S(-3) \stackrel{i} \longrightarrow \operatorname{im}(i) \to 0,
\end{equation}

\begin{equation}\label{sec-seq}
0 \to \operatorname{im}(i) \to  \tilde{Z_{2}} \otimes_A S(-2)  \stackrel{j} \longrightarrow \operatorname{im}(j) \to 0,
\end{equation}

\begin{equation}\label{third-seq}
0 \to \operatorname{im}(j) \to  \tilde{Z_{1}} \otimes_A S(-1) \stackrel{k} \longrightarrow \operatorname{im}(k) \longrightarrow 0,
\end{equation}

\begin{equation}\label{last-seq}
0 \to \operatorname{im}(k) \to   S \longrightarrow \operatorname{Sym}(\langle J,h \rangle) \longrightarrow 0.
\end{equation}
Since $0 \longrightarrow C_{1} \stackrel{\sigma_1} \longrightarrow C_{0}  \to 0$ is a free resolution of $\coker(\sigma_1)$, the mapping cone construction applied to the short exact sequence \eqref{first-seq} yields a graded free resolution of $\operatorname{im}(i)$. In particular, $\operatorname{im}(i)$ admits a graded free resolution of the form
$$
0 \longrightarrow S(-3)  \longrightarrow S(-5) \oplus S(-4) \oplus  S(-3)^7\longrightarrow S(-4)^6 \oplus  S(-3)^{15} \longrightarrow  S(-3)^{19} \longrightarrow 0.
$$
Let a graded free resolution of $\tilde{Z_{2}}= Z_2(\langle J,h \rangle, A)$ be given by
\begin{equation*} \label{res_Z2-2}
\mathcal {R}_{\bullet} : \quad \cdots 
\stackrel{\vartheta_3}{\longrightarrow}  A^{\beta_{3}} 
\stackrel{\vartheta_2}{\longrightarrow}  A^{\beta_{2}}  
\stackrel{\vartheta_1}{\longrightarrow}  A^{\beta_{1}} \longrightarrow 0.
\end{equation*}
Applying the mapping cone construction to the short exact sequence \eqref{sec-seq} then provides a graded free resolution of $\operatorname{im}(j)$, which can be written in the form
\[
\cdots \longrightarrow
\begin{aligned}[t]
	& S(-2)^{\beta_5} \\
	\oplus\;& S(-3)
\end{aligned}
\longrightarrow
\begin{aligned}[t]
	&S(-2)^{\beta_4}\\
	\oplus\;&S(-3)^7\\
	\oplus\;&S(-4)\\
	\oplus\;&S(-5)
\end{aligned}
\longrightarrow
\begin{aligned}[t]
	&S(-2)^{\beta_3}\\
	\oplus\;&S(-3)^{15}\\
	\oplus\;&S(-4)^6
\end{aligned}
\longrightarrow
\begin{aligned}[t]
	&S(-2)^{\beta_2}\\
	\oplus\;&S(-3)^{19}
\end{aligned}
\longrightarrow
S(-2)^{\beta_1}
\longrightarrow 0 .
\]
Since $Z_1(\langle J,h \rangle,A)$ is isomorphic to the first syzygy module of $\langle J,h \rangle$, a truncation of the Koszul complex on $\langle J,h \rangle$, yields a free resolution 
$$
0 \longrightarrow A^4 \longrightarrow A^4 \oplus A^{5} \longrightarrow 0
$$
of $\tilde{Z_{1}}= Z_1(\langle J,h \rangle,A)$. Applying the mapping cone construction to the short exact sequence \eqref{third-seq}, we obtain a graded free resolution of $\operatorname{im}(k)$ given by the complex
\[
\cdots \longrightarrow
\begin{aligned}[t]
& S(-2)^{\beta_5} \\	
\oplus\;& S(-3)
\end{aligned}
\longrightarrow
\begin{aligned}[t]
&S(-2)^{\beta_4}\\
\oplus\;&S(-3)^7\\
\oplus\;&S(-4)\\
\oplus\;&S(-5)
\end{aligned}
\longrightarrow
\begin{aligned}[t]
&S(-2)^{\beta_3}\\
\oplus\;&S(-3)^{15}\\
\oplus\;&S(-4)^6
\end{aligned}
\longrightarrow
\begin{aligned}[t]
&S(-2)^{\beta_2}\\
\oplus\;&S(-3)^{19}
\end{aligned}
\longrightarrow
\begin{aligned}[t]
&S(-2)^{\beta_1} \\ 
\oplus\;&S(-1)^4
\end{aligned}
\longrightarrow 
S(-1)^{9}
\longrightarrow
0.
\]

Finally, using the short exact sequence \eqref{last-seq}, we deduce that the following complex provides an explicit graded free $S$-resolution of the symmetric algebra $\operatorname{Sym}(\langle J,h \rangle)$:

\[
\cdots \longrightarrow
\begin{aligned}[t]
	& S(-2)^{\beta_5} \\
	\oplus\;& S(-3)
\end{aligned}
\longrightarrow
\begin{aligned}[t]
	&S(-2)^{\beta_4}\\
	\oplus\;&S(-3)^7\\
	\oplus\;&S(-4)\\
	\oplus\;&S(-5)
\end{aligned}
\longrightarrow
\begin{aligned}[t]
	&S(-2)^{\beta_3}\\
	\oplus\;&S(-3)^{15}\\
	\oplus\;&S(-4)^6
\end{aligned}
\longrightarrow
\begin{aligned}[t]
	&S(-2)^{\beta_2}\\
	\oplus\;&S(-3)^{19}
\end{aligned}
\longrightarrow
\begin{aligned}[t]
	&S(-2)^{\beta_1} \\ 
	\oplus\;&S(-1)^4
\end{aligned}
\longrightarrow 
S(-1)^{9}
\longrightarrow
S
\longrightarrow
0 .
\]	
\vspace{2mm}

\noindent{\bf Statements and Declarations.} The authors have no competing interests to declare that are relevant to the content of this article.

\vspace{2mm}

\noindent{\bf Acknowledgment.} The first author gratefully acknowledges partial financial support from the Anusandhan National Research Foundation (ANRF), Government of India, via Core Research Grant No. CRG/2023/007668.

\end{document}